\documentclass{article} 
\usepackage[margin=0.5in]{geometry}
\usepackage{amsmath}
\usepackage{amsthm}
\usepackage{amssymb}
\usepackage{authblk}
\usepackage{titling}
\usepackage{hyperref}
\usepackage{indentfirst} 
\usepackage{tikz}
\usetikzlibrary{positioning}
\newtheorem{theorem}{Theorem}
\newtheorem{definition}[theorem]{Definition}
\newtheorem{lemma}[theorem]{Lemma}
\newtheorem{remark}[theorem]{Remark}
\numberwithin{theorem}{section}
\numberwithin{equation}{section}
\title{An Upper Bound for the Menchov-Rademacher Operator for Right Triangles}
\author{Armen Vagharshakyan}
\date{}
\affil{\small\textit{Institute of Mathematics}, \\\small \textit{Armenian National Academy of Sciences}}
\begin{document}
\maketitle
\abstract{
The Menchov-Rademacher inequality is an inequality in harmonic analysis that bounds the $L_2$ norm of a certain maximal operator. It was first established in order to prove almost everywhere convergence of a one-parameter series of orthogonal functions. When two-parameter series of orthogonal functions is considered, the exact way the series is grouped becomes essential. We will consider grouping of a two-parameter series, generated by a sequence of right triangles with a vertex at the origin, who might be non-equilateral, and prove almost everywhere convergence when the eccentricity of those triangles is bounded. In order to carry out the proof, we will derive an analogue of the Menchov-Rademacher inequality for right triangles.
}
{\let\thefootnote\relax\footnote{{Research was supported by the Science Committee of Armenia, grant 18T-1A081.}}}
\section{Introduction}
Certain special orthonormal systems and expansions of functions with respect to them have appeared in the works of D. Bernoulli, F. Bessel, L. Euler, P. Laplace, A. Legendre.
 J. Fourier's method for solving the boundary value problems of mathematical physics, P. Chebyshev's creation of the general theory of orthogonal polynomials (that had arisen from his reasearch on interpolation and the problem of moments), and  D. Hilbert's general theorems on the expansion of functions in a series with respect to an orthonormal system (that had arisen from his research on integral equations) were some of the early incentives to discuss general orthonormal systems of functions. See \cite{EM} for a brief introduction into orthonormal systems.
\newline\indent
Define $d$-parameter systems of orthonormal functions (on a general metric space)  as follows:
\begin{definition}
Denote by $\mathbb{O}^d$ the family of $d$-parameter
systems of functions $\lbrace f_{\vec i}\rbrace_{\vec i\in\mathbb{N}^d}$ that are orthonormal on some metric space $\mathbb{X}$ equipped with a Radon measure $\mu,$ that is:
\begin{align*}
f_{\vec i}\colon \mathbb{X}\rightarrow \mathbb{R}\textit{ for } \vec i\in \mathbb{N}^d,
\\
\int_{\mathbb{X}} f_{ \vec i}f_{\vec i^{\prime}}d\mu=1 \textit{ if } \vec i=\vec i^{\prime}, 
\\
\int_{\mathbb{X}}  f_{\vec i}f_{\vec i^{\prime}}d\mu=0 \textit{ if }  \vec i\neq \vec i^{\prime}.
\end{align*}
\end{definition}
\begin{remark}
The actual choice of the metric space $\mathbb{X}$ equipped with a Radon measure $\mu$ is irrelevant to this paper. The prerequisite that $\mu$ is a Radon measure will be used to ensure the application of B. Levi's theorem in \eqref{e_leviconclusion}.  We refer to \cite{K} for a general discussion of Radon measures. For convenience, the reader may pick $\mathbb{X}=[0,1]^d$ and $\mu$ to be the Lebesgue measure.
\end{remark}
In his study  \cite{Fatou} of trigonometric series, Fatou looked for their convergence criteria expressed in terms of the magnitude of their coefficents, or assuming that the convergence takes place, for properties of functions thus defined. It is in this study that the concept of almost everywhere convergence in relation to trigonometric series first appeared. As it turned out, certain conditons on the magnitude of the coefficients imply almost everywhere convergence of trigonometric series.
\newline\indent
Consider the analogous problem for general orthogonal series, i.e. the problem of finding criteria for almost everywhere convergence for general orthogonal series, expressed  in terms of the magnitude of their coefficents. The following theorem was first proved by Rademacher \cite{Rademacher} (theorem 2 on page 120) and independently by Menchov \cite{Menchov}: 
\begin{theorem}\label{t_mri}
Let $\lbrace f_k\rbrace_{k\in\mathbb{N}}\in\mathbb{O}^1$ be one-parameter orthonormal system of functions. If
\begin{equation}\label{eq_linearmultiplier}
\sum_{k=1}^{+\infty} \left|c_k\right|^2 \ln^2(k+1)<+\infty,    
\end{equation}
then the functional series 
\begin{equation}\label{e_linearseries}
    \sum_{k=1}^{\infty}c_k f_k(x)
\end{equation}
converges almost everywhere on $\mathbb{X}$.
\end{theorem}
\begin{remark} The Menchov-Rademacher theorem \ref{t_mri}  improved an earlier result of Plancherel \cite{Plancherel} (see survey \cite{Ul} for some historical remarks on the Menchov-Rademacher theorem). In turn, Plancherel's result was motivated by questions raised by Cantor and du Bois-Reymond on uniqueness of a function's expansion into orthogonal series (as explained in \cite{Plancherel2}). Theorem \ref{t_mri} can be used to obtain such a uniqueness result. Indeed, assume that the Lebesgue measure of the set $\mathbb{X}$ is finite, so that the $L_1$norms of the functions $f_k$ are uniformly bounded:
\begin{equation*}
\int_{\mathbb{X}}|f_k|d\mu\leq \sqrt{|X|}\left(\int_{\mathbb{X}}|f|^2d\mu\right)^{1/2}= \sqrt{|\mathbb{X}|},\quad k\in\mathbb{N}.
\end{equation*}
Then due to condition \eqref{eq_linearmultiplier}, we can perform term-by-term integration of the series \eqref{e_linearseries}. Consequently, if we denote the sum of the series \eqref{e_linearseries} by $f(x),$
then the coefficients $c_k$ may be restored uniquely from $f$ by the formula
\begin{equation*}
    c_k=\int_{\mathbb{X}}f(x)f_k(x)dx.
\end{equation*}
\end{remark}
Theorem \ref{t_mri} was proved in \cite{Rademacher} and \cite{Menchov} by means of the following lemma:
\begin{lemma}[Menchov-Rademacher inequality for intervals]\label{l_mri}
There exists an absolute constant $\alpha_1>0$ so that for any $\lbrace f_k\rbrace_{k\in\mathbb{N}}\in \mathbb{O}^1$ and any numerical sequence $\lbrace c_k\rbrace_{k\in\mathbb{N}}$ we have
\begin{equation}\label{mrintervals}
    \int_{\mathbb{X}}
    \sup_{1\leq n\leq \mathbb{N}}
    \left|
    \sum_{k=1}^n 
    c_k f_k(x)
    \right|^2d\mu
    \leq  
    \alpha_1^2
    \cdot
    \ln^2 N 
    \cdot
    \left(
    \sum_{k=1}^N
    \left|c_k\right|^2
    \right).
\end{equation}
\end{lemma}
\begin{remark}\label{r_forintervals}
We call inequality \eqref{mrintervals} the Menchov-Rademacher inequality for intervals,
as the summation under the integral sign in \eqref{mrintervals} is carried over natural indices $k$ 
lying in an interval $[1,n].$ In literature, this inequality is known simply as the Menchov-Rademacher inequality, without the word `intervals' appended at the end. The reason we emphasize the word `intervals' is to distinguish this inequality from the one we prove in lemma \ref{l_mrt}. 
\end{remark}
\begin{remark}
In this article, we omit exposing how theorem \ref{t_mri} is derived from lemma \ref{l_mri} in \cite{Rademacher} to avoid redundancy. 
Indeed, in section \ref{s_theoremproof} in a similar, albeit a bit more elaborate manner, theorem \ref{t_mrt} is derived from lemma \ref{l_mrt}. 
\end{remark}
The following notation will be used throughout the article:
\begin{definition}\label{d_tri}
For $0\leq a,b<+\infty,$ denote by $Tri_{a,b}$ the right triangle with vertices $(0,0),(a,0),(b,0)$ (as shown in figure \ref{f1}).
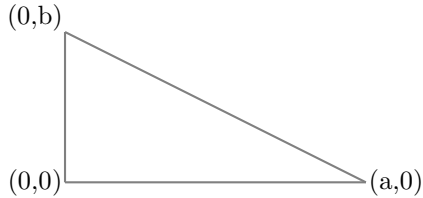
\begin{figure*}[h]
\centering
\begin{tikzpicture}
\draw[gray,thick](0,0) -- (4,0);
\draw[gray,thick](0,0) -- (0,2);
\draw[gray,thick](0,2) -- (4,0);
\node at (4.4,0) {(a,0)};
\node at (-0.4,2.2) {(0,b)};
\node at (-0.4,0) {(0,0)};
\end{tikzpicture}
\caption{triangle $Tri_{a,b}$}
\label{f1}
\end{figure*}
\end{definition}
In this paper, for a two-parameter system of orthonormal functions $\lbrace f_{i,j}\rbrace_{(i,j)\in\mathbb{N}^2}\in\mathbb{O}^2$
and a two-parameter numerical sequence $\lbrace c_{i,j}\rbrace_{(i,j)\in\mathbb{N}^2},$ we are interested in almost everywhere convergence of the functional series 
\begin{equation}\label{eq_functionalseries2d}
    \sum_{(i,j)\in \mathbb{N}^2}c_{i,j}f_{i,j}(x)
\end{equation}
when it is grouped by right triangles $Tri_{a,b}$.
The following theorem regarding almost everywhere convergence of the grouping of the two-parameter series \eqref{eq_functionalseries2d}, generated by a sequence of right triangles $Tri_{a,b}$ with bounded eccentricity, may be viewed as a two-dimensional triangular analogue of theorem \ref{t_mri}:
\begin{theorem}\label{t_mrt}
Let $\lbrace f_{i,j}\rbrace_{(i,j)\in\mathbb{N}^2}\in \mathbb{O}^2$ and let the numerical sequence $\lbrace c_{i,j}\rbrace_{(i,j)\in\mathbb{N}^2}$ satisfy the estimates
\begin{equation}\label{eq_weylmultiplier}
    \sum_{(i,j)\in \mathbb{N}^2}
    \left|(i,j)\right|^{\frac{\ln(2+\sqrt{3})}{\ln(2)}}
    \left|c_{i,j}\right|^2
    <+\infty,
\end{equation}
Consider a sequence of right triangles $\lbrace Tri_{a_k,b_k}\rbrace_{k\in\mathbb{N}}$ (see definition \ref{d_tri}) whose
eccentricity $\gamma$ is bounded, that is
\begin{equation}\label{eq_boundedeccentricity}
    1\leq \sup_{k\in\mathbb{N}}\left(\frac{a_k}{b_k},\frac{b_k}{a_k}\right)=\gamma<+\infty 
\end{equation}
Then the limit
\begin{equation}\label{e_limit}
    \lim_{k\rightarrow +\infty}
    \sum_{(i,j)\in Tri_{a_k,b_k}}\;
    c_{i,j}f_{i,j}(x),
\end{equation}
corresponding to the grouping of the functional series   \eqref{eq_functionalseries2d} by the triangles $Tri_{a_k,b_k},$    
exists and is finite almost everywhere on $\mathbb{X}.$
\end{theorem}
In section \ref{s_theoremproof} we will prove theorem \ref{t_mrt} by means of deriving the following inequality that may be viewed as a triangular analogue of lemma \ref{l_mri}:
\begin{lemma}[Menchov-Rademacher inequality for right triangles]\label{l_mrt}
There exists an absolute constant $\alpha_2>0$ so that for any two-parameter family of orthonormal functions $\lbrace f_{i,j}\rbrace_{(i,j)\in\mathbb{N}^2}\in \mathbb{O}^2$ and any two-parameter numerical sequence $\lbrace c_{i,j}\rbrace_{(i,j)\in\mathbb{N}^2}$ we have
\begin{equation}\label{eq_mrt}
    \int_{\mathbb{X}}
    \sup_{Tri_{a,b}\subset Tri_{N,N}}
    \left|
    \sum_{(i,j)\in Tri_{a,b}} 
    c_{i,j} f_{i,j}(x)
    \right|^2d\mu
    \leq  
    \alpha_2^2
    \cdot
    N^{\frac{\ln{(2+\sqrt{3})}}{\ln{2}}} 
    \cdot
    \left(
    \sum_{(i,j)\in Tri_{N,N}}
    \left|c_{i,j}\right|^2
    \right).
\end{equation}
\end{lemma}
\begin{remark}
We called the inequality \eqref{eq_mrt} the Menchov-Rademacher inequality for right triangles, as in comparison with remark \ref{r_forintervals}, this time the summation is carried over indices $(i,j)$ lying in the right triangle $Tri_{a,b}.$
\end{remark}
\section{Remarks and Complements}
The following two remarks complement the classic Menchov-Rademacher inequality \ref{l_mri}:
\begin{remark}
The best known estimate on the constant $\alpha_1$ appearing in lemma \ref{l_mri} is provided by 
W. Bednorz \cite{Bednorz}:
\begin{equation*}
\alpha_1\leq\frac{1}{3},
\end{equation*}
estimates for $\alpha_1$ obtained earlier being $\alpha_1\leq \ln{2}/\ln{3}$ by E. Kounias \cite{Kounias}, and $\alpha_1\leq 1/(2\ln{2})$ by S. Chobanyan \cite{Chobanyan}. 
\end{remark}
\begin{remark}
The factor $ \ln^2(k+1)$ appearing in condition \eqref{eq_linearmultiplier} is optimal, see the paper \cite{AP} of A. Paszkiewicz for a detailed discussion. 
\end{remark}
The following two remarks are complements to the results \ref{t_mrt} and \ref{l_mrt} that we prove:
\begin{remark}
Unlike the condition \eqref{eq_linearmultiplier} of theorem \ref{t_mri}, the author does not know if the condition \eqref{eq_weylmultiplier} is optimal. 
\end{remark}
\begin{remark}
Our interest in groupings of two-parameter series generated by a sequence of right triangles with a vertex at the origin, who might be non-equilateral, rather than in groupings by different sequences of shapes that cover the set $\mathbb{N}^2$ of all index pairs, is due to the following three observations: 
\\\indent
a) Let $\mathbb{X}=[0,2\pi)^2,$  let $\mu$ be the two-dimensional Lebesgue measure, and $\lbrace f_{i,j}\rbrace$ be the exponential basis $f_{i,j}(x,y)=e^{\textbf{i}(ix+jy)}.$
Assume that
\begin{equation}\label{e_squareconv}
\sum_{(i,j)\in\mathbb{N}^2}\left|c_{i,j}\right|^2<+\infty.
\end{equation} 
Under these assumptions, the result \cite{Fconv} of Ch. Fefferman, is equivalent to claiming that, if  $\lbrace Tri_{a_k,b_k}\rbrace_{k\in\mathbb{N}}$ is a sequence of equilateral right triangles, that is, if the following additional restriction holds:
\begin{equation*}
a_k=b_k,\quad\text{for } k\in\mathbb{N},
\end{equation*}
then the limit \eqref{e_limit} exists almost everywhere. See \cite{Tevzadze},\cite{Sjolin},\cite{Antonov},\cite{MRP} for further refinements of this result.
On the other hand, in a different article \cite{Fdiv}, Ch. Fefferman proved that if no additional restriction is put on the sequence of right triangles $\lbrace Tri_{a_k,b_k}\rbrace_{k\in\mathbb{N}},$
then the condition \eqref{e_squareconv} alone is not sufficient to claim almost everywhere existence of the limit \eqref{e_limit}. Comparing these results with ours, we see that,
on one hand, the assumption \eqref{eq_boundedeccentricity} of bounded eccentricity in our theorem \ref{t_mrt} is intermediate between the assumptions on shapes in these two results of Ch. Fefferman,
on the other hand, our theorem \ref{t_mrt} applies to general orthonormal systems, not just exponentials.
\\\indent
b) Ch. Fefferman's result \cite{Fconv} applies not only to equilateral triangles but also to more general sequences of shapes $\lbrace P_n\rbrace_{n\in \mathbb{N}},$ namely, to those of homothetic copies of polygons. That is, each $P_n$ is a dilation of
the polygon $P_1$ by the same amount in both directions $x$ and $y.$
As for sequences of non-homothetic shapes; the reason to discuss right triangles with a vertex at the origin, rather than other more general sequences of non-homothetic shapes, is the following. Right triangles $Tri_{a,b}$ are dilations of the right triangle $Tri_{1,1}$ in two directions only, namely, each triangle $Tri_{a,b}$ is a dilation of the triangle $Tri_{1,1}$ by amount $a$ in direction $x$ and by amount $b$ in direction $y.$
Consequently, the proofs will minimize difficulties that would arise from comparing those shapes between each other.
\\\indent
c) As for sequences of general triangles; the reason to discuss right triangles with a vertex at the origin, rather than more general sequences of triangles, is the following. Let $T$ be some triangle in the first quadrant who is not necessarily a right triangle with a vertex at the origin. Consider a sequence of triangles $\lbrace T_n\rbrace_{n\in\mathbb{N}}$ generated by dilations of the triangle $T$ by amount $a_n>0$ in direction $x$ and by amount $b_n>0$ in direction $y.$ Assume that for the sequence of triangles $\lbrace T_n\rbrace_{n\in\mathbb{N}}$ assumption \eqref{eq_boundedeccentricity} of theorem \ref{t_mrt} holds.  
Then for the sequence of triangles $\lbrace T_n\rbrace_{n\in\mathbb{N}}$ to cover the first quadrant, and thus to be a valid 
choice for generating a grouping of a two-parameter series, it is necessary that $\lbrace T_n\rbrace_{n\in\mathbb{N}}$ is a sequence of right triangles with a vertex at the origin. 
\end{remark}
The following two remarks single out some of the techniques used in proving lemma \ref{l_mrt}:
\begin{remark}\label{r_recursive}
The following theorem in the Analysis of Algorithms elicits
asymptotic estimates from recursive estimates:
\begin{theorem}[Master Theorem in the Analysis of Algorithms \cite{Bentleyetal}]\label{t_master}
Let a function $T\colon N\rightarrow \mathbb{R}$ satisfy the bound
\begin{equation*}
T(n)\leq aT(n/b)+f(n),
\end{equation*}
and let
\begin{equation*}
f(n)=O\left(n^c\right).
\end{equation*}
If $c<\log_b (a),$ then
\begin{equation*}
T(n)=O\left(n^{\log_b(a)}\right).
\end{equation*}
If $c>\log_b (a),$ then
\begin{equation*}
T(n)=O\left(n^{c}\right).
\end{equation*}
\end{theorem}
\end{remark}
\begin{remark}\label{r_rquadratic}
It is a well-known fact of Linear Algebra that the maximum of a quadratic form on the unit sphere equals its largest eigenvalue.
\end{remark}
\section{Proof of lemma \ref{l_mrt}}
\begin{figure}[h]
\centering
\begin{tikzpicture}
\draw[gray,thick](0,0) -- (6,0);
\draw[gray,thick](0,0) -- (0,4);
\draw[gray,thick](0,4) -- (6,4);
\draw[gray,thick](6,0) -- (6,4);
\draw[gray,thick](2,0) -- (2,4);
\draw[gray,thick](2,4) -- (5,0);
\draw[gray,thick](2,0) -- (5,0);
\draw[gray,thick](2,2) -- (3.5,2);
\draw[gray,thick](3.5,0) -- (3.5,2);
\node at (2.5,2.5) {$T^1$};
\node at (4,0.5) {$T^2$};
\node at (2.5,0.5) {$T^3$};
\node at (6.2,0) {m};
\node at (0,4.2) {n};
\end{tikzpicture}
\caption{triangles $T^1,T^2,T^3$} 
\label{f3}
\end{figure}
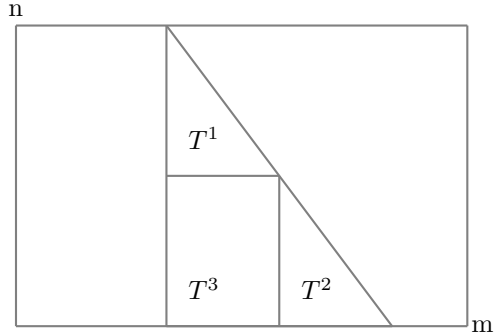

\noindent
Let $S$ be a family of sets in $\mathbb{R}^d.$ Introduce the generalized Menchov-Rademacher operator corresponding to the family of sets $S,$ as a maximal operator acting on a family of orthogonal functions $\lbrace f_{\vec i}\rbrace_{\vec i\in \mathbb{N}^d}$ and on a numerical sequence $\lbrace c_{\vec i}\rbrace_{\vec i\in \mathbb{N}^d}$ as follows:
\begin{equation*}
\left(
\lbrace f_{\vec i}\rbrace_{\vec i\in \mathbb{N}^d}
,
\lbrace c_{\vec i}\rbrace_{\vec i\in \mathbb{N}^d}
\right)
\rightarrow \sup_{I\in S}\left(\sum_{\vec i\in I\cap \mathbb{N}^d}c_{\vec i}f_{\vec i}(x)\right).
\end{equation*}
Denote by $mr(S)$ the sharp upper bound for the $L_2$ norm of the generalized Menchov-Rademacher operator, that is
\begin{equation*}
mr(S)=
\sup_{\lbrace f_{\vec i}\rbrace\in \mathbb{O}^d}\;\sup_{\sum \left|c_{\vec i}\right|^2 \leq 1}\;
\left|\left| \sup_{I\in S}\left|\sum_{\vec i\in I\cap \mathbb{N}^d}c_{\vec i}f_{\vec i}(x)\right|\right|\right|_{L_2\left(\mathbb{X},\mu\right)}.
\end{equation*}
For brevity, we introduce the following notations:\\
a) $SINT_n$ is the family of intervals lying in $[1,n]$:
$
SINT_n=\lbrace [k_1,k_2] \colon k_1,k_2\in[1,n]\rbrace,
$
\\
b) $INT_n$ is a subset of $SINT_n$ and consists of those intervals lying in $[1,n]$ who have $1$ as an endpoint:
$
INT_n=\lbrace [1,k] \colon k\in [1,n]\rbrace,
$
\\
c) $TRI_n$ is the family of right triangles with a vertex at the origin, lying in the triangle $Tri_{n,n}$:
$
TRI_n=\lbrace Tri_{a,b} \colon Tri_{a,b}\subset Tri_{n,n} \rbrace,
$
\\
d) $HREC_{n,m}$ is a family of rectangles of fixed height:
$
HREC_{n,m}=\lbrace [x,y]\times [0,m]\colon x,y\in [0,n]\rbrace
,$
\\
e) $HTRI_{n,m}$ is a family of triangles of fixed height:
$
HTRI_{n,m}=\lbrace a+T_{b,m}\colon 0\leq a\leq a+b\leq n\rbrace.
$\\
In terms of these notations, lemma \ref{l_mri} may be paraphrased as:
\begin{equation}\label{eq_para}
    mr(INT_n)=O(\ln(n)),
\end{equation}
and the claim of lemma \ref{l_mrt} may be paraphrased as:
\begin{equation}\label{l_concise}
mr(TRI_{n})=O\left(n^{\frac{ln(2+\sqrt{3})}{2\ln(2)}}\right).
\end{equation}
We split the proof of the estimate \eqref{l_concise} into three steps.
Namely, we first obtain an estimate for $mr(HREC_{m,n}),$ corresponding to rectangles of fixed height; then obtain an estimate for $mr(HTRI_{m,n}),$ 
corresponding to triangles of fixed height; and finally derive an estimate for $mr(TRI_n),$ corresponding to right triangles.
\\
\textbf{Step 1.}
By \eqref{eq_para} and the triangle inequality, we have
\begin{equation}\label{step11}
    mr(SINT_n)\leq 2 \cdot mr(INT_n),\quad \textit{for }n\in\mathbb{N}.
\end{equation}
Now let $\lbrace f_{i,j} \rbrace_{(i,j)\in\mathbb{N}^2}\in \mathbb{O}^2$ and
\begin{equation}\label{cond_c}
\sum_{(i,j)\in  \mathbb{N}^2}\left|c_{i,j}\right|^2\leq 1.
\end{equation}
For the family of orthogonal functions
\begin{equation}\label{step12}
g_i=\sum_{j=1}^m c_{i,j}f_{i,j},\quad n\in\mathbb{N},
\end{equation}
using inequality \eqref{cond_c}, we have
\begin{equation}\label{step13}
\left|\left|\sup_{I\in SINT_n}\left(\sum_{i\in I} g_i\right)\right|\right|_{L_2(\mathbb{X},\mu)}
\leq 
mr\left(SINT_n\right),\quad \textit{for }n\in\mathbb{N}.
\end{equation}
By combining \eqref{step12} and \eqref{step13} we get
\begin{equation}\label{e:est_REC}
mr(HREC_{n,m})\leq  mr(SINT_n),\quad \text{for }n,m\in\mathbb{N}.
\end{equation}
\textbf{Step 2.}
First, we will obtain the recursive estimate \eqref{e:htri_rec} for $mr(HTRI_{m,n}).$
The following identity will serve as the initial condition for the recursive estimate:
\begin{equation}\label{e:htri_HREC_init}
mr(HTRI_{n,1})=mr(SINT_n),\quad \text{for }n\in\mathbb{N}.
\end{equation}
Now let $\lbrace f_{i,j} \rbrace\in \mathbb{O}^2$ be a two-parameter system of orthonormal functions, and let the numerical sequence $\lbrace c_{i,j}\rbrace$ satisfy the condition
\begin{equation}\label{c_coeffs}
\sum_{i=1}^n\sum_{j=1}^m \left|c_{i,j}\right|^2\leq 1.
\end{equation}
Split the sum \eqref{c_coeffs} into two parts (the lower part and the upper part) as follows:
\begin{align}\label{d_p}
p^2_1=
\sum_{i=1}^n\sum_{j=1}^{\lfloor m/2\rfloor}
\left|c_{i,j}\right|^2,\quad
\textit{taking } p_1>0,
\\
p^2_2=
\sum_{i=1}^n\sum_{j=\lfloor m/2\rfloor+1}^{m}
\left|c_{i,j}\right|^2,\quad
\textit{taking }
p_2>0,
\end{align}
so that the condition \eqref{c_coeffs} translates into
\begin{equation}\label{e:circle1}
p_1^2+p_2^2\leq 1.
\end{equation}
Split a triangle $T\in HTRI_{n,m}$ into three parts - $T^1,T^2,T^3$ as shown in figure \ref{f3}; to be precise:
\begin{align*}
T^1=\lbrace (x,y)\in T\colon       y>m/2 \rbrace,
\\
T^2=\lbrace (x,y)\in T\colon    x>(a+b)/2,  y\leq m/2 \rbrace,
\\
T^3=\lbrace (x,y)\in T\colon     a\leq x\leq (a+b)/2,  y\leq m/2 \rbrace.
\end{align*}
\newline
The following inclusions hold:
\begin{align*}
T^1-(a,m/2)\subset HTRI_{n,m/2},
\\
T^2-((a+b)/2,0)\subset HTRI_{n,m/2},
\\
T^3-(a,0)\subset HREC_{n,m/2}.
\end{align*}
Using these inclusions and definitions \eqref{d_p}, by the triangle inequality we have:
\begin{equation*} 
\left|\left|\sup_{T\in HTRI_{n,m} }\left|\sum_{(i,j)\in T}c_{i,j}f_{i,j}\right|\right|\right|_{L_2(\mathbb{X},\mu)}
\leq 
p_2 \cdot mr\left(HTRI_{n,\frac{m}{2}}\right)
+
p_1 \cdot mr\left(HREC_{n,\frac{m}{2}}\right)
+
p_1 \cdot mr\left(HTRI_{n,\frac{m}{2}}\right)
\end{equation*}
Squaring this inequality and maximizing it over condition \eqref{e:circle1}, we get
\begin{gather}\label{e_recursivee}
mr(HTRI_{n,m})^2\leq 
\left(mr\left(HTRI_{n,\frac{m}{2}}\right)+mr\left(HREC_{n,\frac{m}{2}}\right)\right)^2+mr^2\left(HTRI_{n,\frac{m}{2}}\right).
\end{gather}
Consequently, the following recursive linear estimate, that is weaker than the non-linear estimate \eqref{e_recursivee}, holds:
\begin{equation}\label{e:htri_rec}
mr\left(HTRI_{n,m}\right)\leq \sqrt{2}\cdot mr\left(HTRI_{n,\frac{m}{2}}\right)+mr\left(HREC_{n,\frac{m}{2}}\right)
\end{equation}
We now apply the estimate \eqref{e:htri_rec} recursively over $m,$ starting with the initial condition \eqref{e:htri_HREC_init}, to get
\begin{equation*}
mr\left(HTRI_{n,2^k}\right)
\leq 
2^{k/2}\cdot mr(SINT_n)
+
\sum_{l=1}^{k-1}2^{(k-1-l)/2}\cdot mr\left(HREC_{n,2^l}\right)
\end{equation*}
Now, we insert the estimate \eqref{e:est_REC} obtained in Step $1,$ to get
\begin{equation*}
mr\left(HTRI_{n,2^k}\right)
\leq  
\left(2^{k/2} +\sum_{l=1}^{k-1}2^{(k-1-l)/2}\right)
\cdot 
mr(SINT_n)
\quad \text{for }n,k\in\mathbb{N}.
\end{equation*}
Further, by \eqref{eq_para} and \eqref{step11},
\begin{equation}\label{e_left}
mr\left(HTRI_{n,2^k}\right)
\leq
O\left(2^{k/2}
\cdot 
\ln(n)
\right)
,
\quad \text{for }n,k\in\mathbb{N}.
\end{equation}
Apparently,
\begin{equation*}
HTRI_{n,m}\subset HTRI_{n,2^{\lceil \log_2(m)\rceil}}
\end{equation*}
so that 
\begin{equation}\label{e_subset}
mr(HTRI_{n,m})\leq mr\left(HTRI_{n,2^{\lceil \log_2(m)\rceil}}\right). 
\end{equation}
We use inequality \eqref{e_subset} to substitute the power of two - $2^k$ on the left side of inequality \eqref{e_left} to an arbitrary natural number $m.$ This way we obtain  the following bound for triangles of fixed height:
\begin{equation}\label{e:est_htri} 
mr(HTRI_{n,m})=O(\sqrt{m}\ln (n)),\quad \quad\text{for }n,m\in\mathbb{N}.
\end{equation}
\textbf{Step 3.}
Let $\lbrace f_{i,j} \rbrace\in \mathbb{O}^2$ be a two-parameter system of orthonormal functions, and let numerical coefficients $\lbrace c_{i,j}\rbrace$ satisfy the condition
\begin{equation}\label{e_restriction}
\sum_{i,j\in Tri_{n,n}} \left|c_{i,j}\right|^2\leq 1.
\end{equation}
Split the triangle $Tri_{n,n}$ into four triangles $T^1,T^2,T^3,T^4$ as shown in figure \ref{f4}; to be precise:
\begin{align*}
T^1=Tri_{\frac{n}{2},\frac{n}{2}}+\left(0,\frac{n}{2}\right),
\\
T^2=Tri_{\frac{n}{2},\frac{n}{2}}+\left(\frac{n}{2},0\right),
\\
T^3=Tri_{\frac{n}{2},\frac{n}{2}},
\\
T^4=Tri_{n,n}\setminus\left(T^1\cup T^2 \cup T^3\right).
\end{align*}
\begin{figure}[h]
\centering
\begin{tikzpicture}
\draw[gray,thick](0,0) -- (6,0);
\draw[gray,thick](0,0) -- (0,4);
\draw[gray,thick](3,0) -- (3,2);
\draw[gray,thick](0,2) -- (3,2);
\draw[gray,thick](0,4) -- (6,0);
\draw[gray,thick](0,2) -- (3,0);
\node at (0.5,0.5) {$T^3$};
\node at (3.5,0.5) {$T^2$};
\node at (0.5,2.5) {$T^1$};
\node at (2.5,1.5) {$T^4$};
\end{tikzpicture}
\caption{triangles $T^1,T^2,T^3,T^4$} \label{f4}
\end{figure}
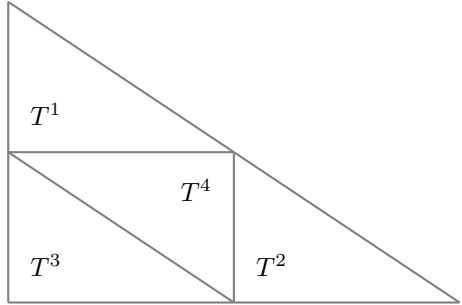
\begin{equation*}
p_k^2=\sum_{(i,j)\in T^k\cap \mathbb{N}^2}\left|c_{i,j}\right|^2,\quad \textit{for } k=1,2,3,4, \quad \textit{taking } p_k\geq 0.
\end{equation*}
Condition \eqref{e_restriction} translates into
\begin{equation}\label{e:circle2}
\sum_{k=1}^4 p_k^2\leq 1.
\end{equation}
\begin{remark}\label{r_shape}
Note that a triangle $T=T_{a,b}\in TRI_n,$ depending on whether $a>n/2$ or $a\leq n/2$ and whether
$b>n/2$ or $b\leq n/2,$ may be of either of the following four forms:
\\
a) $T$ is a triangle in $T^3.$ 
\\
b) $T$ is a union of a triangle in $T^1,$ a rectangle in $HREC_{\frac{n}{2},\frac{n}{2}},$ and a triangle in $HTRI_{\frac{n}{2},\frac{n}{2}},$ 
\\
c) $T$ is a symmetry w.r.t the line $y=x$ of a set described in the case 2 above,
\\
d) $T$ is a union of a triangle in $T^1,$ a triangle in $T^2,$ and of a set-theoretical difference of two shapes - the square $[0,\frac{n}{2}]\times[0,\frac{n}{2}]$ and a triangle in $(\frac{n}{2},\frac{n}{2})-TRI_{\frac{n}{2},\frac{n}{2}}.$ 
\end{remark}
\begin{remark}
As it will be seen from the estimate \eqref{e_est} of quadratic form, case d) of remark \ref{r_shape} is the case that provides the largest contribution to the exponent
\begin{equation*}
\frac{\ln{(2+\sqrt{3})}}{\ln{2}}
\end{equation*}
that appears in lemma \ref{l_mri}. This is because in case d) the triangle $T$ splits into the most number of parts. 
\end{remark}
As a consequence of remark \ref{r_shape} that allows four forms for the triangle $T$, by the triangle inequality we have the following estimate that has four corresponding terms on the right hand side:
\begin{gather*}
\left|\left|\sup_{T\in TRI_{n} }\left|\sum_{(i,j)\in T}c_{i,j}f_{i,j}\right|\right|\right|_{L_2(\mathbb{X},\mu)}^2
\leq 
\\
\leq
p_3^2 \cdot
mr^2\left(TRI_{\frac{n}{2}}\right)
+
\left(
p_1 \cdot
mr\left(TRI_{\frac{n}{2}}\right)
+\sqrt{p_3^2+p_4^2} \left(
mr\left(HREC_{\frac{n}{2},\frac{n}{2}}\right)
+
mr\left(HTRI_{\frac{n}{2},\frac{n}{2}}\right)\right)
\right)^2
+
\\
+
\left(
p_2 
\cdot
mr\left(TRI_{\frac{n}{2}}\right)+\sqrt{p_3^2+p_4^2} \left(mr\left(HREC_{\frac{n}{2},\frac{n}{2}}\right)
+mr\left(HTRI_{\frac{n}{2},\frac{n}{2}}\right)\right)
\right)^2+
\left((p_1+p_2+p_4) \cdot mr\left(TRI_{\frac{n}{2}}\right)+\sqrt{p_3^2+p_4^2}\right)^2
\end{gather*}
If we rearrange this estimate by powers of $mr(TRI_n)$ and apply the estimates \eqref{e:est_REC} and \eqref{e:est_htri} obtained 
for $mr(HREC_{n,m})$ and $mr(HTRI_{n,m})$ in steps $1$ and $2$ of this section, then we get:
\begin{equation}\label{e_nearlythere}
mr(TRI_{n})^2\leq 
\left(p_1^2+p_2^2+p_3^2+(p_1+p_2+p_4)^2\right)\cdot mr^2\left(TRI_{\frac{n}{2}}\right)
+\alpha_3 \sqrt{n}\ln(n)\cdot mr\left(TRI_{\frac{n}{2}}\right)+\alpha_4 n\ln^2(n)
\end{equation}
for some constants $0<\alpha_3,\alpha_4<+\infty$.
\begin{remark}\label{r_quadratic}
For a moment, let us discuss the factor
\begin{equation}\label{e_factor}
    p_1^2+p_2^2+p_3^2+(p_1+p_2+p_4)^2
\end{equation}
that appears on the right hand side of \eqref{e_nearlythere}.
We can interpret the inequality \eqref{e:circle2} as a claim that the vector $(p_1,p_2,p_3,p_4)$ belongs to the unit ball. This is why, by remark \ref{r_rquadratic}, we can estimate 
the quadratic form \eqref{e_factor}
of the variables $p_1,p_2,p_3,p_4$ by its largest eigenvalue $\kappa=2+\sqrt{3}$ as follows: 
\begin{equation}\label{e_est}
p_1^2+p_2^2+p_3^2+(p_1+p_2+p_4)^2\leq \kappa=2+\sqrt{3}.
\end{equation}
\end{remark}
Insert estimate \eqref{e_est} of remark \ref{r_quadratic} into inequality
\eqref{e_nearlythere}
to obtain the following linear recursive estimate:
\begin{equation*} 
mr(TRI_{n})\leq  
\sqrt{\kappa}\cdot mr\left(TRI_{\frac{n}{2}}\right)+\alpha_5\sqrt{n}\ln(n),\quad \textit{for some constant }0<\alpha_5<+\infty.
\end{equation*}
The proof of estimate \eqref{l_concise} (and consequently of lemma \ref{l_mrt}) now follows from theorem \ref{t_master}.
\section{Proof of theorem \ref{t_mrt}}\label{s_theoremproof}
\noindent
We split the proof of theorem \ref{t_mrt} into three steps. 
Namely, we first use condition \eqref{eq_weylmultiplier} to prove existence of limit almost everywhere when the grouping of the series \eqref{eq_functionalseries2d} is by a sequence of right rectangles whose endpoints grow with geometric speed. Then we use lemma \ref{l_mrt} to compare that grouping with the grouping introduced in theorem \ref{t_mrt}. And finally we use B. Levi's theorem to prove almost everywhere convergence of theorem \ref{t_mrt}.\\
\textbf{Step 1.}
Due to condition \eqref{eq_weylmultiplier} on the numerical sequence $\lbrace c_{i,j}\rbrace,$ the assumptions of theorem \ref{t_mri} are satisfied for the functional series
\begin{equation}\label{e_s_all}
\sum_{k=1}^{\infty}\left(\sum_{(i,j)\in Tri_{k+1,k+1}\setminus Tri_{k,k}}
c_{i,j}f_{i,j}(x)
\right).
\end{equation}
Hence, the limit of the partial sums of the series \eqref{e_s_all},
\begin{equation*}
  \lim_{k\rightarrow +\infty}
    \sum_{(i,j)\in Tri_{k,k}}\;
 c_{i,j}f_{i,j}(x)
\end{equation*}
exists and is finite for almost every $x\in \mathbb{X}.$
In particular, using $\gamma$ for eccentricity (see \eqref{eq_boundedeccentricity}), the following limit over a subsequence of $\lbrace Tri_{k,k}\rbrace_{k=1}^{+\infty},$
\begin{equation}\label{e_gammas}
   \lim_{p\rightarrow +\infty}
    \sum_{(i,j)\in Tri_{\gamma^p,\gamma^p}}\;
 c_{i,j}f_{i,j}(x)
\end{equation}
exists and is finite for almost every $x\in \mathbb{X}.$
\\
\textbf{Step 2.}
In this step we prove inequality \eqref{e_tail} for  $p\in\mathbb{N}.$ Colloquially, 
inequality \eqref{e_tail} estimates the difference between a member of the sequence appearing in the limit \eqref{e_limit} and
the corresponding member of \eqref{e_gammas}.\\
Note that
by the bound on eccentricity \eqref{eq_boundedeccentricity}, we have
\begin{equation}\label{e_inclusion}
\left[
\gamma^p\leq min(a_k,b_k)<\gamma^{p+1}
\right]
\implies
\left[
Tri_{a_k,b_k}
\subset 
Tri_{\gamma^{p+2},\gamma^{p+2}}
\right].
\end{equation} 
Also note that for a pair of indices $(i,j)\in \mathbb{N}^2,$ from geometric considerations we have
\begin{equation}\label{e_linetonorm}
    \left[(i,j)\not\in Tri_{\gamma^p,\gamma^p}\right]
\implies 
\left[|(i,j)|\geq \gamma^p/\sqrt{2}\right].
\end{equation}
and

\begin{equation}\label{e_linetonorm2}
    \left[(i,j)\in Tri_{\gamma^{p+2},\gamma^{p+2}}\right]
\implies 
\left[|(i,j)|\leq \gamma^{p+2}\right].
\end{equation}
Denote by $\lbrace c^{\ast}_{i,j}\rbrace_{(i,j)\in Tri_{\gamma^{p+2},\gamma^{p+2}}}$ the restriction of the numerical sequence $\lbrace c_{i,j}\rbrace$ to the set of indices $Tri_{\gamma^{p+2},\gamma^{p+2}}\setminus Tri_{\gamma^{p},\gamma^{p}},$ that is
\begin{align}\label{e_stars}
c^{\ast}_{i,j}=c_{i,j}& \quad \textit{for }(i,j)\in Tri_{\gamma^{p+2},\gamma^{p+2}}\setminus Tri_{\gamma^{p},\gamma^{p}},
\\
c^{\ast}_{i,j}=0& \quad \textit{for }(i,j)\in Tri_{\gamma^{p},\gamma^{p}}.\nonumber
\end{align}
We have
\begin{align}\label{e_tail}
    \int_{\mathbb{X}}
    \max_{k\colon \gamma^p\leq min(a_k,b_k)<\gamma^{p+1}}
    \left|
    \sum_{(i,j)\in Tri_{a_{k},b_{k}}\setminus Tri_{\gamma^p,\gamma^p}}
    \;
    c_{i,j}f_{i,j}(x)
    \right|^2
d\mu
=
    \int_{\mathbb{X}}
    \max_{k\colon \gamma^p\leq min(a_k,b_k)<\gamma^{p+1}}
    \left|
    \sum_{(i,j)\in Tri_{a_{k},b_{k}}}
    \;
    c^{\ast}_{i,j}f_{i,j}(x)
    \right|^2
d\mu
\leq\nonumber
\\
\leq
  \int_{\mathbb{X}}
    \max_{T\in TRI_{\gamma^{p+2},\gamma^{p+2}}}
    \left|
    \sum_{(i,j)\in T}
    \;
    c^{\ast}_{i,j}f_{i,j}(x)
    \right|^2 d\mu
\leq\nonumber
\\
\leq
\alpha_6
 \left(
\sum_{(i,j)\in Tri_{\gamma^{p+2},\gamma^{p+2}}}
|c^{\ast}_{i,j}|^2
\right)
\cdot
\left(\gamma^{p+2}\right)^{\frac{\ln(2+\sqrt{3})}{\ln(2)}}
\leq
\alpha_7 \gamma^{\frac{2\ln(2+\sqrt{3})}{\ln(2)}}\cdot
\sum_{(i,j)\in Tri_{\gamma^{p+2},\gamma^{p+2}}\setminus Tri_{\gamma^p,\gamma^p}}
\left(
|c_{i,j}|^2
|(i,j)|^{\frac{\ln(2+\sqrt{3})}{\ln(2)}}
\right)
\end{align}
for some absolute constants $0<\alpha_6,\alpha_7<+\infty.$ 
Indeed, the first equality in \eqref{e_tail} follows from \eqref{e_stars} and \eqref{e_inclusion},                           
the next inequality in \eqref{e_tail} follows from \eqref{e_inclusion},
the next inequality in \ref{e_tail} follows from lemma \ref{l_mrt},
and the last inequality in \ref{e_tail} follows from \eqref{e_stars} and \eqref{e_linetonorm}.
\\
\textbf{Step 3.} We summate both sides of the inequality \eqref{e_tail} over $p\in\mathbb{N}$ in order to get                        
\begin{align}\label{e_beforelevi}
   \sum_{p=1}^{+\infty} \int_{\mathbb{X}}
    \max_{k\colon \gamma^p\leq min(a_k,b_k)<\gamma^{p+1}}
    \left|
    \sum_{(i,j)\in Tri_{a_{k},b_{k}}\setminus Tri_{\gamma^p,\gamma^p}}
    \;
    c_{i,j}f_{i,j}(x)
    \right|^2
d\mu
    \leq\nonumber
\\
\leq
\alpha_7 \gamma^{\frac{2\ln(2+\sqrt{3})}{\ln(2)}}\cdot
\sum_{(i,j)\in \mathbb{N}^2}
\left(
|c_{i,j}|^2
|(i,j)|^{\frac{\ln(2+\sqrt{3})}{\ln(2)}}
\mathbf{card}
\left[
p\in\mathbb{N}\colon (i,j)\in Tri_{\gamma^{p+2},\gamma^{p+2}}\setminus Tri_{\gamma^p,\gamma^p} 
\right]
\right)
\end{align}
where $\textbf{card}$ denotes the cardinality. 
Now note that due to \eqref{e_linetonorm} and\eqref{e_linetonorm2}, the function
\begin{equation*}
(i,j)\rightarrow \mathbf{card}
\left[
p\in\mathbb{N}\colon (i,j)\in Tri_{\gamma^{p+2},\gamma^{p+2}}\setminus Tri_{\gamma^p,\gamma^p} 
\right]
\end{equation*}
is bounded by a constant that depends on $\gamma$ only. Using this fact and the assumption \eqref{eq_weylmultiplier} of theorem \ref{t_mrt}, inequality \eqref{e_beforelevi} implies:
\begin{equation}\label{e_leviassumption}
   \sum_{p=1}^{+\infty} \int_{\mathbb{X}}
    \max_{k\colon \gamma^p\leq min(a_k,b_k)<\gamma^{p+1}}
    \left|
    \sum_{(i,j)\in Tri_{a_{k},b_{k}}\setminus Tri_{\gamma^p,\gamma^p}}
    \;
    c_{i,j}f_{i,j}(x)
    \right|^2
d\mu<+\infty
\end{equation}
By B. Levi's theorem applied to Radon measure $\mu,$ the condition \eqref{e_leviassumption} implies that
\begin{equation}\label{e_leviconclusion}
    \max_{k\colon \gamma^p\leq min(a_k,b_k)<\gamma^{p+1}}
    \left|
    \sum_{(i,j)\in Tri_{a_{k},b_{k}}\setminus Tri_{\gamma^p,\gamma^p}}
    \;
    c_{i,j}f_{i,j}(x)
    \right|\rightarrow 0,\quad\text{for a.e. }x\in\mathbb{X}.
\end{equation}
Finally, \eqref{e_gammas} and \eqref{e_leviconclusion} imply the conclusion \eqref{e_limit} of theorem \ref{t_mrt}.

\end{document}